\documentclass{amsart}

\usepackage{amsmath,amssymb,amsfonts,enumerate,amsthm,graphicx,color}

\newcommand{\q}{\mbox{q}\,}

\newcommand{\E}{\mbox{E}}

\newcommand{\V}{\mbox{V}}

\newtheorem{thm}{Theorem}[section]
\newtheorem{cor}[thm]{Corollary}
\newtheorem{lem}[thm]{Lemma}
\newtheorem{prop}[thm]{Proposition}

\numberwithin{equation}{section}

\begin{document}
\bibliographystyle{amsplain}

%date{}
\title{Arithmetical rank of the cyclic and bicyclic graphs }

\author{Margherita Barile}

\address{Margherita Barile\\ Dipartimento di Matematica, Universit\`{a} di Bari ``Aldo Moro'', Via E.
Orabona 4, 70125 Bari, Italy. }

\email{barile@dm.uniba.it}

\author{Dariush Kiani}

\address{Dariush Kiani\\Department of Mathematics, Amirkabir University of technology, Tehran, Iran\\
and School of Mathematics, Institute for Research in Fundamental
Sciences (IPM).}

\email{dkiani@aut.ac.ir}

\author{Fatemeh Mohammadi}

\address{Fatemeh Mohammadi\\Department of Mathematics, Amirkabir University of technology, Tehran, Iran.}

\email{f-mohammadi@cic.aut.ac.ir}

\author{Siamak Yassemi}
\address{Siamak Yassemi\\Department of Mathematics, University of
Tehran, Tehran, Iran and School of Mathematics, Institute for
Research in Fundamental Sciences (IPM).}

\email{yassemi@ipm.ir}

%\thanks{Emails: barile@dm.uniba.it, dkiani@aut.ac.ir, f-mohammadi@cic.aut.ac.ir, yassemi@ipm.ir}
\thanks{D. Kiani was supported in part by a grant from IPM (No. 87200116)}
\thanks{S. Yassemi was supported in part by a grant from
IPM (No. 87130211).}

\keywords{Arithmetical rank, projective dimension, edge ideals,
set-theoretic complete intersection ideals}

\subjclass[2000]{13A15; 13F55; 13D02; 05C38} \maketitle
\begin{abstract}

\noindent We show that for the edge ideals of the graphs consisting
of one cycle or two cycles of any length connected through a vertex
or a path, the arithmetical rank equals the projective dimension.

\end{abstract}

%%%%%%%%%%%%%%%%%%%%%%%%%%%%%%%%%%%%%%%%%%%%%%%%%%%%%%%%%%%%%%%%%%%%%%%%%%%%%%%%%%%%%%%%%%%%%%%%%%%%%%%%%%%%%%%%%%%%%%
\section{Introduction}

For any homogeneous ideal $I$ of a polynomial ring
$R=K[x_1,\ldots,x_n]$ there exists a  \textit{graded} minimal finite
free resolution $$ 0\rightarrow \bigoplus_j
R(-j)^{\beta_{pj}}\rightarrow \ldots\rightarrow \bigoplus_j
R(-j)^{\beta_{1j}}\rightarrow R\rightarrow R/I\rightarrow 0$$ of
$R/I$, in which $R(-j)$ denotes the graded free module obtained by
shifting the degrees of elements in $R$ by $j$. The numbers
$\beta_{ij}$, which we shall refer to as the $i$th Betti numbers of
degree $j$ of $R/I$, are independent of the choice of the graded minimal
finite free resolution. We also define the $i$th Betti number of $I$
as $\beta_i:=\sum \beta_{ij}$.

Given a polynomial ring $R$ over a field, and a graph $G$ having the
set of indeterminates as its vertex set $V(G)$ and the set of edges
$E(G)$, one can associate with $G$ a monomial ideal of $R$: this
ideal is generated by the products of the vertices of each edge in
$E(G)$, and is hence generated by squarefree quadratic monomials. It
is called the {\it edge ideal} $I(G)$ of $G$, and has been
intensively studied by Simis, Vasconcelos and Villarreal in
\cite{SVV}.
 The arithmetical rank ($ara$), i.e., the least number of elements of
 $R$ which generate a given monomial ideal up to radical,
 is in general bounded below by its projective dimension ($pd$), i.e.,
 by the length of every minimal free resolution of the quotient of $R$
  with respect to the ideal.
 The simplicial complex
$\Delta_G$ of a graph $G$ is defined by
$$\Delta_G=\{A\subseteq V(G)|A \,\, \mbox{is an independent set in}\,\, G\},$$
where $A$ is an independent set in $G$ if none of its elements are
adjacent. Note that $\Delta_G$ is precisely the Stanley-Reisner
simplicial
complex of $I(G)$. \\
For any simplicial complex $\Delta$ on the vertex set $V(\Delta)$,
the \textit{Alexander dual} of $\Delta$ is the simplicial complex
defined by
$${\Delta}^{\ast}:=\{F\subseteq V(\Delta)|\ V(\Delta)\setminus F\notin \Delta\}.$$ The link of a face $F\in \Delta$ is defined as the simplicial
complex $$Link_{\Delta}F:=\{G\in \Delta | G\cup F\in \Delta\ and\
G\cap F=\emptyset\}.$$
  In recent times, the projective dimension has been determined for
  large classes of edge ideals, where it is independent of the ground field:
  in Jacques' thesis it was computed for acyclic graphs (see also \cite{Jac2}),
  but also for the graphs $C_n$, consisting of one cycle of length
  $n$. Jacques, in \cite[Theorem
6.1.8]{Jac1}, using Hochster's formula~\cite{Ho}, showed that
for a graph $G$, the Betti numbers are \\
\\
(*) $\ \ \ \ \ \ \ \ \ \ \ \ \ \ \ \ \ \beta_{i,d}(G)=\sum_{H\subset
G,|V(H)|=d} dim_k \tilde{H}_{i-2}(\varepsilon (H);K).$
\\
\\ Then he used (*) for providing formulas for the graded Betti
numbers of special classes of graphs including lines, cycles and
complete graphs. He proved the following theorems.

\noindent{\bf Theorem A} \cite[Lemma 8.2.7]{Jac1} The reduced
homology of the disjoint union of the cyclic graph
    $C_n$ and any non empty graph $G$ may be expressed as follows:\[
{\tilde{H}}_i( \varepsilon (C_n\cup G);K)= \left\lbrace
  \begin{array}{c l}
    {\tilde{H}}_{i-\frac{2n+1}{3}}( \varepsilon (G);K), & \text{if $n\equiv 1\ {\rm mod}\ 3$}\\
{\tilde{H}}_{i-\frac{2n-1}{3}}( \varepsilon (G);K), & \text{if
$n\equiv 2\ {\rm mod}\ 3$}.
  \end{array}
\right. \]

\noindent{\bf Theorem B} \cite[Corollary 7.6.30]{Jac1} The non zero
Betti numbers in degree $n$ and the projective dimension of $C_n$ in
degree $n$ are the following:
\\
\\
$\beta_{\frac{2n}{3},n}=2$,  \ \ \ and\  $pd\ I(C_n)=\frac{2n}{3}$,
if
$n\equiv 0\ {\rm mod}\ 3$, \\
 $\beta_{\frac{2n+1}{3},n}=1$,\ \  and\  $pd\ I(C_n)=\frac{2n+1}{3}$, if $n\equiv
1\ {\rm mod}\ 3$,
\\
$\beta_{\frac{2n-1}{3},n}=1$,\ \  and\  $pd\ I(C_n)=\frac{2n-1}{3}$,
if $n\equiv 2\ {\rm mod}\ 3$.
\\

\noindent{\bf Theorem C} \cite[Lemma 8.1.3]{Jac1} The reduced
homology of the disjoint union of the line graph
    $L_n$ and any non empty graph $G$ may be expressed as follows:\[
{\tilde{H}}_i( \varepsilon (L_n\cup G);K)= \left\lbrace
  \begin{array}{c l}
    {\tilde{H}}_{i-\frac{2n}{3}}( \varepsilon (G);K), & \text{if $n\equiv 0\ {\rm mod}\ 3$}\\
0, & \text{if $n\equiv 1\ {\rm mod}\ 3$}\\
{\tilde{H}}_{i-\frac{2n-1}{3}}( \varepsilon (G);K), & \text{if
$n\equiv 2\ {\rm mod}\ 3$}.
  \end{array}
\right. \] From the proof of \cite[Corollary 7.7.35]{Jac1} one can derive the following result.
\\

\noindent{\bf Theorem D} \cite[Corollary 7.7.35]{Jac1} The
projective dimension of the line graph is independent of the
characteristic of the chosen field and is \[ pd\ I(L_n) =
\left\lbrace
  \begin{array}{c l}
    \frac{2n}{3} & \text{if $n\equiv 0\ {\rm mod}\ 3$}\\
\frac{2n-2}{3} & \text{if $n\equiv 1\ {\rm mod}\ 3$}\\
\frac{2n-1}{3} & \text{if $n\equiv 2\ {\rm mod}\ 3$}.
  \end{array}
\right. \]
 All Betti numbers of
$L_n$ in degree $n$ are zero if $n \equiv 1\ {\rm mod}\ 3$. Otherwise
the non zero Betti numbers of degree $n$ of $L_n$
are $$\beta_{\frac{2n}{3} ,n}\ I(L_n) = 1, \rm{if}\ \ {\it n} \equiv 0\
{\rm mod}\ 3,$$
$$\beta_{\frac{2n-1}{3} ,n} \ I(L_n)= 1, \rm{if}\ \ {\it n} \equiv 2\ {\rm mod}\ 3.$$ In \cite{Cor}
  an explicit formula is given for the Betti numbers of a special kind of bipartite graphs, the so-called {\it Ferrers graphs}.
  In \cite{B2} it is shown that the arithmetical rank equals the projective dimension for a special class
  of acyclic graphs,
   in \cite{B4} that this is also true for all Ferrers graphs.
   In the present paper we  prove that the same equality holds for all cyclic and bicyclic
   graphs. By \textit{bicyclic
   graph} we mean a graph which consists of two cycles that have exactly one vertex in common or are connected by a path.
 In particular, we will see that the projective dimension of the edge ideals of these graphs
 does not depend on the characteristic of the ground field.
   %%%%%%%%%%%%%%%%%%%%%%%%%%%%%%%%%%%%%%%%%%%%%%%%%%%%%%%%%%%%%%%%%%%%%%%%%%%%%%%%%%%
\section{The arithmetical rank of cyclic graphs}
Let $K$ be a field, and consider the polynomial ring $R=K[x_1,\dots,
x_n]$, where $n\geq3$. Let $C_n$ be the graph on the vertex set
$\{x_1, \dots, x_n\}$ whose set of edges is $\{\{x_1, x_2\}, \{x_2,
x_3\}, \dots, \{x_{n-1}, x_n\}, \{x_1, x_n\}\}$. Then its edge ideal
is the following ideal of $R$:
$$I(C_n)=(x_1x_2, x_2x_3, \dots, x_{n-1}x_n, x_1x_n).$$
We will show that for all $n$, $pd\ I(C_n)=\, ara\,I(C_n)$. In
general, for any ideal $I$ of $R$ we have that $cd\ I\leq\,ara\,I$,
where $cd$ denotes the local cohomological dimension (see \cite{H},
Example 2, p.~414) and, whenever $I$ is a
monomial ideal, $pd\ I=\, cd\,I$ (see \cite{L1}, Theorem 1). Hence
it will suffice to show that, for all $n$, $ara\ I(C_n)\leq pd\
I(C_n)$, i.e., to produce $ pd\ I(C_n)$ elements of $R$ generating
$I(C_n)$, up to radical. Among the available tools, we have, on the
one hand, Jacques' result providing explicit formulas for the
projective dimension of $I(C_n)$.

On the other hand, we know that a finite set of elements of $R$
which generate a given ideal up to radical can be constructed
according to the following criterion, which is due to Schmitt and
Vogel.
\begin{lem}\label{lemma}{\rm (\cite{Sch}, p.~249)} Let $P$ be a finite subset of elements of $R$. Let $P_0,\dots, P_r$ be subsets of $P$ such that
\begin{list}{}{}
\item[(i)] $\bigcup_{i=0}^rP_i=P$;
\item[(ii)] $P_0$ has exactly one element;
\item[(iii)] if $p$ and $p'$ are different elements of $P_i$ $(0<i\leq r)$ there is an integer $i'$ with $0\leq i'<i$ and an element in $P_{i'}$ which divides $pp'$.
\end{list}
\noindent We set $q_i=\sum_{p\in P_i}p^{e(p)}$, where $e(p)\geq1$
are arbitrary integers. We will write $(P)$ for the ideal of $R$
generated by the elements of $P$.  Then we get
$$\sqrt{(P)}=\sqrt{(q_0,\dots,q_r)}.$$
\end{lem}
\par\medskip\noindent
We have to distinguish between three cases, depending on the residue
of $n$ modulo 3. The cases $n\equiv 0,1$ mod 3 can be settled by a
direct application of Lemma \ref{lemma}. The case $n\equiv 2$ mod 3
is more interesting, since it needs some additional non trivial
computations on the generators.
\begin{prop}\label{0} Suppose that $n=3m$, for some integer $m$. Set $q_0=x_1x_2$, $q_1=x_1x_{3m}+x_2x_3$, and,
for $1\leq i\leq m-1$, set
\begin{eqnarray*} q_{2i}&=&x_{3i+1}x_{3i+2}\\
q_{2i+1}&=&x_{3i}x_{3i+1}+x_{3i+2}x_{3i+3}.
\end{eqnarray*}
Then
$$I(C_n)=\sqrt{(q_0,\dots, q_{2m-1})}.$$
In particular, $ara\,I(C_n)\leq 2m.$
\end{prop}
\begin{proof}
For all $i=0,\ldots, m-1$, the monomial $q_{2i}$ divides the product
of the two summands of $q_{2i+1}$. By Lemma \ref{lemma} it follows
that
$$(x_{3i}x_{3i+1}, x_{3i+1}x_{3i+2}, x_{3i+2}x_{3i+3})=\sqrt{(q_{2i}, q_{2i+1})}.$$
 This implies the claim.
\end{proof}
Using the same arguments as in the proof of Proposition~\ref{0},
from Lemma \ref{lemma} we can deduce the next result.
\begin{prop}\label{1} Suppose that $n=3m+1$, for some integer $m$. Set $q_0=x_1x_2$, $q_1=x_1x_{3m+1}+x_2x_3$, and,
for $1\leq i\leq m-1$, set
\begin{eqnarray*} q_{2i}&=&x_{3i+1}x_{3i+2}\\
q_{2i+1}&=&x_{3i}x_{3i+1}+x_{3i+2}x_{3i+3},
\end{eqnarray*}
and, finally, $q_{2m}=x_{3m}x_{3m+1}$. Then
$$I(C_n)=\sqrt{(q_0,\dots, q_{2m})}.$$
In particular, $ara\,I(C_n)\leq 2m+1.$
\end{prop}

\begin{prop}\label{2} Suppose that $n=3m+2$, for some integer $m$. Set $q_0=x_1x_2$, $q_1=x_2x_3+x_4x_5$, and, for
$1\leq i\leq m-1$, set
\begin{eqnarray*} q_{2i}&=&x_{3i}x_{3i+1}+x_{3i+2}x_{3i+3}\\
q_{2i+1}&=&x_{3i+2}x_{3i+3}+x_{3i+4}x_{3i+5},
\end{eqnarray*}
and, finally, $q_{2m}=x_1x_{3m+2}+x_{3m}x_{3m+1}$. Then
$$I(C_n)=\sqrt{(q_0,\dots, q_{2m})}.$$
In particular, $ara\,I(C_n)\leq 2m+1.$
\end{prop}
\begin{proof}
The claim for $m=1$ was proven in \cite{B2}, Example 1. So let
$m\geq 2$. Set $J_m=(q_0, \dots, q_{2m})$. It suffices to show that
$I(C_n)\subset\sqrt{J_m}$. In  this proof, for all $f,g\in R$, by
abuse of notation we will write $f\equiv g$ whenever $f-g$ or $f+g$
belongs to $J_m$, and,   $f \equiv_{q_i} g$ whenever $f-g$ or $f+g$
is divisible by $q_i$. In this way, $f\equiv g$ or $f \equiv_{q_i}
g$  assures that $f\in J_m$ occurs if and only if $g\in J_m$. We
first show that
\begin{equation}\label{eq0} x_1^{2^m}x_{3m+2}^{2^{m+1}}\in J_m.\end{equation}
Set
\begin{eqnarray*}
u_m&=&x_1^{2^{m-1}}x_{3m+2}^{2^m},\\
v_m&=&x_3x_4x_5\prod_{i=2}^mx_{3i}^{3\cdot 2^{i-2}},\\
w_m&=&(x_{3m}x_{3m+1}x_{3m+2})^{2^{m-1}}.
\end{eqnarray*}
We prove that
\begin{equation}\label{eq1} u_m\equiv v_m\equiv w_m.\end{equation}
Note that $x_1x_{3m+2}v_m$ is a multiple of $x_1x_{3m+2}x_4x_5$, and
$$x_1x_{3m+2}x_4x_5\equiv_{q_0}x_1x_{3m+2}(x_2x_3+x_4x_5)\in J_m,$$
whence we deduce that $x_1x_{3m+2}v_m\in J_m$. Thus (\ref{0}) will
imply that
$$x_1^{2^m}x_{3m+2}^{2^{m+1}}=x_1x_{3m+2}u_m\in J_m,$$
as claimed in (\ref{eq0}). We prove (\ref{eq1}) by induction on
$m\geq 2$. First take $m=2$. We have $q_2=x_3x_4+x_5x_6$,
$q_3=x_5x_6+x_7x_8$ and $q_4=x_1x_8+x_6x_7$, so that
$$
v_2=x_3x_4x_5x_6^3\equiv_{q_2}x_5^2x_6^2x_6^2\equiv_{q_3}x_6^2x_7^2x_8^2=w_2\\
\equiv_{q_4}x_1^2x_8^4=u_2,$$ which shows (\ref{eq1}) for $m=2$. Now
suppose that $m>2$ and that the claim is true for $m-1$. We have:
\begin{eqnarray*} v_m&=&v_{m-1}x_{3m}^{3\cdot 2^{m-2}}\equiv w_{m-1}x_{3m}^{3\cdot 2^{m-2}}\\\\
&=&(x_{3m-3}x_{3m-2}x_{3m-1})^{2^{m-2}}x_{3m}^{3\cdot 2^{m-2}}
=(x_{3m-3}x_{3m-2})^{2^{m-2}}x_{3m-1}^{2^{m-2}}x_{3m}^{3\cdot 2^{m-2}}\\\\
&\equiv_{q_{2m-2}}&(x_{3m-1}x_{3m})^{2^{m-2}}x_{3m-1}^{2^{m-2}}x_{3m}^{3\cdot
2^{m-2}}
=(x_{3m-1}x_{3m})^{2^{m-1}}x_{3m}^{2\cdot 2^{m-2}}\\\\
&\equiv_{q_{2m-1}}&(x_{3m+1}x_{3m+2})^{2^{m-1}}x_{3m}^{2^{m-1}}
=(x_{3m}x_{3m+1})^{2^{m-1}}x_{3m+2}^{2^{m-1}}=w_m\\\\
&\equiv_{q_{2m}}&(x_1x_{3m+2})^{2^{m-1}}x_{3m+2}^{2^{m-1}}
=x_1^{2^{m-1}}x_{3m+2}^{2^m}=u_m.
\end{eqnarray*}
This completes the proof of (\ref{eq1}) and of (\ref{eq0}). We have
thus shown that
\begin{equation}\label{eq2} x_1x_{3m+2}\in\sqrt{J_m}.\end{equation}
But then
\begin{equation}\label{eq3} x_{3m}x_{3m+1}=q_{2m}-x_1x_{3m+2}\in\sqrt{J_m}.\end{equation}
In general, whenever, for some $i\in\{2,\dots, m\}$,
\begin{equation}\label{eq4} x_{3i}x_{3i+1}\in\sqrt{J_m},\end{equation}
from the fact that $x_{3i}x_{3i+1}$ divides $x_{3i-1}x_{3i}\cdot
x_{3i+1}x_{3i+2}$, i.e., the product of the summands of $q_{2i-1}$,
by Lemma \ref{lemma} one deduces that
\begin{equation}\label{eq5} x_{3i-1}x_{3i}\in\sqrt{J_m}.\end{equation}
Since $x_{3i-3}x_{3i-2}=q_{2i-2}-x_{3i-1}x_{3i}$, this in turn
implies that
\begin{equation}\label{eq6} x_{3i-3}x_{3i-2}\in\sqrt{J_m}.\end{equation}
Finally, since $x_{3i-3}x_{3i-2}$ divides $x_{3i-4}x_{3i-3}\cdot
x_{3i-2}x_{3i-1}$, i.e., the product of the summands of $q_{2i-3}$,
by Lemma \ref{lemma} we again conclude that
\begin{equation}\label{eq7} x_{3i-2}x_{3i-1}\in\sqrt{J_m}.\end{equation}
Therefore, since (\ref{eq4}) implies (\ref{eq5}), (\ref{eq6}) and
(\ref{eq7}), for all $i=2,\dots, m$, from (\ref{eq3}) one can derive
by descending induction on $h$, that $x_hx_{h+1}\in \sqrt{J_m}$ for
all $h=3,\dots, 3m+1$. In particular we have that
$x_3x_4\in\sqrt{J_m}$, which, together with $q_1\in J_m$, yields
$x_2x_3\in\sqrt{J_m}$ by Lemma \ref{lemma}. This, together with
(\ref{2}) and $q_0\in J_m$, shows that $I(C_n)\subset\sqrt{J_m}$, as
claimed.
\end{proof}
Theorem B and Propositions~\ref{0}, \ref{1}, \ref{2} imply our
main result.
\begin{thm}\label{theorem} Let $n\geq 3$ be an integer. Then
$$ara\ I(C_n)=\, pd\ I(C_n)=\left\{\begin{array}{cll}
\frac{2n}3&{\rm if}&n\equiv 0\ {\rm mod}\ 3\\\\
\frac{2n+1}3&{\rm if}&n\equiv 1\ {\rm mod}\ 3\\\\
\frac{2n-1}3&{\rm if}&n\equiv 2\ {\rm mod}\ 3.
\end{array}
\right.$$
\end{thm}
\par\medskip\noindent
Every ideal $I(C_n)$ is of pure height $\lceil \frac{n}2\rceil$,
where $\lceil a\rceil$ denotes the least integer not less than $a$.
Recall that an ideal is called a {\it set-theoretic complete
intersection} if its arithmetical rank equals its height. In view of
Theorem \ref{theorem} we thus have the following.
\begin{cor} $I(C_n)$ is a set-theoretic complete intersection only for $n=3$ and $n=5$.
\end{cor}

%%%%%%%%%%%%%%%%%%%%%%%%%%%%%%%%%%%%%%%%%%%%%%%%%%%%%%%%%%%%%%%%%%%%%%%%%%%%%%%%%%%%%%%%%%%%%%%%%%%%%%%%%%%%%%%%%%%%%%
\section{The arithmetical rank of bicyclic graphs}

In this section by $\equiv$, we mean $\equiv \ ($mod$\
   3)$ and all equivalence relations will be considered modulo $3$.
Let $a_1,\ldots,a_s$ be subsets of the finite set $V$. Define
$\varepsilon(a_1,\ldots,a_s;V)$ to be the simplicial complex which
has vertex set $\bigcup_{i=1}^s (V\setminus a_i)$ and maximal faces
$V\setminus a_1,\ldots,V\setminus a_s$. Let $\Delta=\Delta_G$, and let $F\in {\Delta}^{\ast}$
and $e_1,\ldots,e_r$ be all the edges of $G$ which are disjoint from
$F$. Then
$Link_{{\Delta}^{\ast}}F=\varepsilon(e_1,\ldots,e_r;V(G)\setminus
F)$ by \cite[Proposition 3.3]{Jac2}. According to ~\cite[Proposition
6.1.6]{Jac1},
 associating $F$ with the induced subgraph $H$
of $G$ on the vertex set $V(G)\setminus F$ defines a bijection
between the faces of ${\Delta}^{\ast}$ and the set of induced
subgraphs of $G$ which have at least one edge. Let $H$ be an induced
subgraph of the graph $G$. If $H$ is associated with the face $F$ of
${\Delta}^{\ast}$ as described above, we write $\varepsilon(H)$
for $\varepsilon(e_1,\ldots,e_s;V)$, where $e_1,\ldots,e_s$ are the
edges of $H$ and $V$ is the vertex set $V(G)\setminus F$ (or
equivalently the vertex set of $H$). In this section, using (*), we
find explicit descriptions of the projective dimension of all
bicyclic graphs. For every vertex $u$ of a graph $G$ we denote by
$N_G(u)$ the set of vertices adjacent to $u$. In the proof of our
main results we will use the Mayer-Vietoris sequence for the reduced
homology of simplicial complexes, which, for any pair $\Delta_1,\
\Delta_2$ of simplicial complexes, has the following form (see
\cite[Remark 6.2.13]{Jac1}):
$$\ldots\rightarrow {\tilde{H}}_{i}(\Delta_1 \cap \Delta_2;K)\rightarrow
{\tilde{H}}_{i}(\Delta_1 ;K)\oplus {\tilde{H}}_{i}(
\Delta_2;K)\rightarrow {\tilde{H}}_{i}(\Delta_1 \cup
\Delta_2;K)\rightarrow$$
$$\qquad\qquad\qquad\qquad\qquad\qquad\qquad\qquad\qquad{\tilde{H}}_{i-1}(\Delta_1 \cap
\Delta_2;K)\rightarrow\ldots.$$
\begin{lem}\label{deg one}
For a graph $G$ with an edge $\{u,v\}$ such that $deg(v)=1$, we have
${\tilde{H}}_{i}(\varepsilon(G);K)=0$, if some vertex in $N_G(u)$
has an adjacent vertex of degree one in $G$. Otherwise,
${\tilde{H}}_{i}(\varepsilon(G);K)={\tilde{H}}_{i-t}(\varepsilon(H);K)$,
where $t=|N_G(u)|$ and $H$ is the induced subgraph on
$V(G)\setminus(\{u\}\cup N_G(u))$, provided $H$ is non empty.
\end{lem}
\begin{proof}{In this and in the following proofs we will omit the coefficient
field in the homology groups.
We set $V=V(G)$. Let $N_G(u)=\{v,u_1,\ldots,u_{t-1}\}$ and
$\{u,v\}, \{u,u_1\},\ldots,\{u,u_{t-1}\},e_1,\ldots,e_r$ be the
edges of $G$. We can write $\varepsilon(G)=E_1\cup E_2$, where
$E_1=\varepsilon(\{u,u_1\},\ldots,\{u,u_{t-1}\},e_1,\ldots,e_r;V)$
and $E_2=\varepsilon(\{u,v\};V)$.\\
The intersection of these simplicial complexes is:
\\
$E_1\cap
E_2=\varepsilon(\{u,v,u_1\},\ldots,\{u,v,u_{t-1}\},\{u,v\}\cup
e_1,\ldots,\{u,v\}\cup e_r;V)$
\\
$=\varepsilon(\{u_1\},\ldots,\{u_{t-1}\}, e_1,\ldots, e_r;V\setminus(\{u,v\}))$ (see \cite[Lemma 3.4]{Jac2}).\\
If there exists a vertex $v_i$ of degree one such that $\{u_i,v_i\}\in
E(G)$, then without loss of generality we can assume that $e_1=\{u_i,
v_i\}$. Then $E_1\cap E_2=\varepsilon(\{u_1\},\ldots,\{u_{t-1}\},\\
\{u_i,v_i\}, e_2,\ldots,
e_r;V\setminus(\{u,v\}))=\varepsilon(\{u_1\},\ldots,\{u_{t-1}\},
e_2,\ldots, e_r;V\setminus(\{u,v\}))$, whose reduced homology is
identically zero, since $v_i\in V\setminus(\{u,v\})$ and $v_i$
belongs to all faces of $E_1\cap E_2$. Otherwise,
 by \cite[Lemma 3.5]{Jac2} we have
 ${\tilde{H}}_i(E_1\cap E_2)=\tilde{H}_{i-t+1}(\varepsilon(H))$, for all $i$, where $H$ is the induced
 subgraph on $V\setminus(\{u\}\cup N_G(u))$. Since $E_2$ is a
simplex, ${\tilde{H}}_i( E_2)=0$ for all $i$. Also, ${\tilde{H}}_i(
E_1)=0$ for all $i$, since $v$ belongs to all faces of $E_1$.
Using the Mayer-Vietoris sequence (for $\Delta_i=E_i$) we deduce that
${\tilde{H}}_{i}(\varepsilon(G))={\tilde{H}}_{i-1}(E_1\cap E_2)$,
which completes the proof.
 }\end{proof}

The next result can be deduced from Lemma~\ref{deg one} by a trivial
inductive argument.
\begin{cor}\label{empty}{
Let $n\equiv 0$. Suppose that $L_n$ intersects graph $G$ only at one
of its endpoints. Then, for all $i$, we have
${\tilde{H}}_i(\varepsilon(G\cup L_n))=
{\tilde{H}}_{i-\frac{2n}{3}}(\varepsilon(G\setminus L_n)).$
}\end{cor}
\begin{thm}\label{pd vertex}
Let $G$ be the graph which is a joint of two cycles $C_n$ and $C_m$
in a common vertex. Then the following hold:
\\
\\
$(a)$ If $|V(G)|\equiv 1$, then $pd\ I(G)= ara\ I(G)=\frac{2|V(G)|+1}{3}$.\\
$(b)$ If $|V(G)|\equiv 0$, then $pd\ I(G)= ara\ I(G)=\frac{2|V(G)|}{3}$.\\
$(c)$ If $|V(G)|\equiv 2$, then $pd\ I(G)= ara\
I(G)=\frac{2|V(G)|+2}{3}$, for $m\equiv 0$, whereas $pd\ I(G)= ara\
I(G)=\frac{2|V(G)|-1}{3}$ otherwise.

\end{thm}
\begin{proof}{We will prove the claim by
showing that the desired number is, on the one hand, a lower bound
for $pd\ I(G)$, on the other hand, an upper bound for $ara \ I(G)$.

Let $V=V(G)$. Consider the labeling for $V$ such that
$V(C_n)=\{y_1,y_2,\ldots,y_n\}$, and $V(C_m)=\{x_1,x_2,\ldots,x_m\}$,
where $x_1=y_1$.   Up to exchanging $m$ and $n$ we have the following cases.
\\

\textbf{Case $1$.} $|V|\equiv 0$ or $1$.  \\
First let
$n=3$. Then $m\equiv 1$ or $m\equiv 2$. In view of (*) the $i$th
Betti number of degree $|V|$ is $\beta_{i,|V|}(G)=
 dim_k {\tilde{H}}_{i-2}(\varepsilon(G);K)$.
So we compute the reduced homology of $G$ of degree $|V|$. We can
write
$\varepsilon(G)=\varepsilon(\{x_1,x_2\},\ldots,\{x_m,x_1\},\{x_1,y_2\},\{y_2,y_3\},\{y_{3},x_1\};\\V)=E_1\cup
E_2$, where
$E_1=\varepsilon(\{x_1,x_2\},\ldots,\{x_m,x_1\},\{x_1,y_2\},\{y_{3},x_1\};V)$
and
 $E_2=\varepsilon(\{y_2,y_3\};V)$.\\ By \cite[Lemma 3.4]{Jac2}, the intersection of these simplicial
complexes is:
\begin{eqnarray*} E_1\cap
E_2&=&\varepsilon(\{x_1,x_2,y_2,y_3\},\ldots,\{x_1,x_m,y_2,y_3\},\{x_1,y_2,y_3\};V)\\
&=&\varepsilon(\{x_1\},\{x_2,x_3\},\ldots,\{x_{m-1},x_m\};V\setminus\{y_2,y_3\}).
\end{eqnarray*}
 By
\cite[Lemma 3.5]{Jac2}, ${\tilde{H}}_i(E_1\cap
E_2)=\tilde{H}_{i-1}(\varepsilon(L_{m-1}))$, for all $i$. Since $E_2$ is a
simplex, ${\tilde{H}}_i( E_2)=0$ for all $i$. Applying
Lemma~\ref{deg one} to $E_1$ for $v=y_2$ (and $u=x_1$, so that
$N(u)=\{x_2,x_m,y_2,y_3\}$), we obtain
${\tilde{H}}_i(E_1)=\tilde{H}_{i-4}(\varepsilon(L_{m-3}))$, for all $i$.\\
If $m\equiv 1$, then ${\tilde{H}}_i(E_1)=0$ for all $i$ by Theorem D
and (*) (since $m-3\equiv 1$). By the Mayer-Vietoris sequence we deduce that, for all $i$,
${\tilde{H}}_i(\varepsilon(G))={\tilde{H}}_{i-1}(E_1\cap E_2)$.
Theorem D and (*) then show that $\tilde H_{i-2}(\varepsilon(L_{m-1}))\ne0$ (i.e.,
${\tilde{H}}_i(\varepsilon(G))\neq 0$)
if and only if $i-2+2=\frac{2(m-1)}3$, if and only if $i=\frac{2|V|}{3}-2$. In view of (*) we deduce that
$pd\ I(G)\geq \frac{2|V|}{3}$.
\\
If $m\equiv 2$, then ${\tilde{H}}_i(E_1\cap E_2)=0$ for all $i$ by
Theorem D and (*) (since $m-1\equiv 1$). By the Mayer-Vietoris sequence we deduce that
${\tilde{H}}_i(\varepsilon(G))={\tilde{H}}_{i}(E_1)$ for all $i$. Theorem D and
(*) show that ${\tilde{H}}_i(\varepsilon(G))\neq 0$ if and only if
$i=\frac{2|V|+1}{3}-2$. In view of (*) we deduce that $pd\ I(G)\geq
\frac{2|V|+1}{3}$.
\\
So we can assume that $n\geq 4$. Moreover, since $n$ and $m$ cannot be both divisible by 3, we may assume that $m\equiv 1$ or $m\equiv 2$. In view of (*) we compute the reduced homology of $G$ of degree $|V|$. We can
write\\
$\varepsilon(G)=\varepsilon(\{x_1,x_2\},\ldots,\{x_m,x_1\},\{x_1,y_2\},\{y_2,y_3\},
\ldots,\{y_{n-1},y_n\},\{y_{n},x_1\};V)=\\E_1\cup E_2$, where
$E_1=\varepsilon(\{x_2,x_3\},\ldots,\{x_m,x_1\},\{x_1,y_2\},\ldots,\{y_{n-1},y_n\},\{y_{n},x_1\};V)$
and
 $E_2=\varepsilon(\{x_1,x_2\};V)$. We have that $E_1=\varepsilon(L_m\cup C_n)$, where $L_m: x_2x_3 \ldots x_m x_1$.
\\ The intersection of these simplicial
complexes is $E_1\cap
E_2=$$\varepsilon(\{x_1,x_2,x_3\},\{x_1,x_2,x_3,\\x_4\},\ldots,\{x_1,x_2,x_{m-1},x_m\},\{x_m,x_1,x_2\},\{x_1,y_2,x_2\},
\ldots,\{y_{n-1},y_{n},x_1,x_2\},\{y_{n},x_1,\\x_2\};V)$
$=\varepsilon(\{x_3\},\{x_m\},\{y_2\},\{y_n\},\{x_4,x_5\},\ldots,
\{x_{m-2},x_{m-1}\},\{y_3,y_4\},\ldots,\{y_{n-2},\\y_{n-1}\};V\setminus\{x_1,x_2\})$
(see \cite[Lemma 3.4]{Jac2}). \\
 By
\cite[Lemma 3.5]{Jac2}, ${\tilde{H}}_i(E_1\cap
E_2)=\tilde{H}_{i-4}(\varepsilon(H))$ for all $i$, where $H$ is the
induced subgraph on $V\setminus\{x_1,x_2,x_3,x_{m},y_2,y_n\}$. We
have ${\tilde{H}}_i(E_1\cap
E_2)=\tilde{H}_{i-4}(\varepsilon(L_{m-4}\cup L_{n-3}))$ for all $i$.
Since $E_2$ is a simplex, ${\tilde{H}}_i( E_2)=0$ for all $i$.\\

\textbf{Case $1.1$.} Let $m\equiv 2$.
\\
By Theorem C, ${\tilde{H}}_i(E_1\cap E_2)=0$ for any $i$, since
$m-4\equiv 1$.  Applying Corollary
~\ref{empty} to the path $L_{m-2}:\ x_2x_3\ldots x_{m-1}$, we get
that, for all $i$,
${\tilde{H}}_{i}(E_1)=\tilde{H}_{i-\frac{2(m-2)}{3}}(\varepsilon(L_{2}\cup
C_n ))$, where $L_2:x_1 x_m$. If we apply Lemma~\ref{deg one} once
again for $v=x_m$ (and $u=x_1$, so that $N(u)=\{x_m,y_2,y_n\}$), we
then obtain
${\tilde{H}}_{i}(E_1)=\tilde{H}_{i-\frac{2(m-2)}{3}-3}(\varepsilon(L_{n-3}))$, for all $i$,
where $L_{n-3}:y_3\ldots y_{n-1}$. In part $(a)$, we have $n\equiv
0$. Theorem D and (*) show that ${\tilde{H}}_{i}(E_1)\neq 0$ if and
only if $i=\frac{2|V|+1}{3}-2$. The Mayer-Vietoris sequence  implies that ${\tilde{H}}_i(\varepsilon(G))\neq 0$
if and only if $i=\frac{2|V|+1}{3}-2$. By (*) it follows that $pd\
I(G)\geq\frac{2|V|+1}{3}$, as claimed. In part $(b)$, we have
$n\equiv 2$. Theorem D and (*) show that ${\tilde{H}}_{i}(E_1)\neq
0$ if and only if $i=\frac{2|V|}{3}-2$. As above, it follows that
${\tilde{H}}_i(\varepsilon(G))\neq 0$ if and only if
$i=\frac{2|V|}{3}-2$. In view of (*) we deduce that $pd\
I(G)\geq\frac{2|V|}{3}$, as claimed.
\\

\textbf{Case $1.2$.} Let $m\equiv 1$.
\\
By Theorem C, ${\tilde{H}}_i(E_1\cap
E_2)=\tilde{H}_{i-4-\frac{2(m-4)}{3}}(\varepsilon(L_{n-3}))$ for any
$i$. Moreover, applying Corollary~\ref{empty} to the path
$L_{m-1}:x_2x_3\ldots x_m$, we get that, for all $i$,
${\tilde{H}}_{i}(E_1)=\tilde{H}_{i-\frac{2(m-1)}{3}}(\varepsilon(C_{n}))$.
In part $(a)$, we have $n\equiv 1$. Hence, by Theorem C,
${\tilde{H}}_i(E_1\cap E_2)=0$ for any $i$, since $n-3\equiv 1$. On
the other hand, Theorem B and (*) show that
${\tilde{H}}_{i}(E_1)\neq 0$ if and only if $i=\frac{2|V|+1}{3}-2$.
We deduce that ${\tilde{H}}_{i}(\varepsilon(G))\neq 0$ if and only
if $i=\frac{2|V|+1}{3}-2$. By (*) it follows that $pd \ I(G)\geq
\frac{2|V|+1}{3}$, as claimed. In part $(b)$, we have $n\equiv 0$.
By Theorem D and (*) it follows that ${\tilde{H}}_i(E_1\cap E_2)\neq
0$ if and only if $i=\frac{2|V|}{3}-2$, in which case the $i$th
homology module is equal to $K$. Theorem B and (*) show that
${\tilde{H}}_{i}(E_1)\neq 0$ if and only if $i=\frac{2|V|}{3}-2$ in
which case it is equal to $K^2$. Thus the Mayer-Vietoris sequence
implies that ${\tilde{H}}_{i}(\varepsilon(G))\neq 0$ for
$i=\frac{2|V|}{3}-2$. In view of
(*) we deduce that $pd \ I(G)\geq \frac{2|V|}{3}$, as claimed.
\\

\textbf{Case $2$.} Let $|V|\equiv 2$.\\
\textbf{Case $2.1$.} $m\equiv 0$.
\\
We have $n\equiv 0$.  First assume that $n=m=3$.
 We first
show that in this case $pd\ I(G) \geq 4$. We use the fact that $pd\
I(G) = cd\ I(G)$, (see \cite[Theorem 1]{L2}), where $cd$ denotes the
local cohomological dimension, i.e., for any ideal $I$ of R, $cd\  I$
is the maximum index $i$ for which the local cohomology module $H_I
^{i} (R)$ (of $R$ with respect to $I$) does not vanish. We have that
$I(G) = I \cap J$, where $I = (x_1, x_2x_3, y_2y_3)$ and $J = (x_2,
x_3, y_2, y_3)$. It is well-known that, whenever an ideal is a
complete intersection, its height is the only index for which the
cohomology module of $R$ with respect to this ideal does not vanish (see \cite[Proposition 3.8]{HT} together with \cite[Theorem 4.4]{HT}, or together with \cite[Example 2, p.~414]{H}).
 Since $I + J = (x_1, x_2, x_3, y_2, y_3)$ we thus have that
$H_{I+J} ^{i} (R) \neq 0$ if and only if $i = 5$. We also have that $cd\ J =4$. In the
Mayer-Vietoris sequence for local cohomology (see \cite{HT}, Section
$3$)
$$\ldots\rightarrow H_{I+J} ^{4}(R)\rightarrow H_{I} ^{4} (R)\oplus H_{J} ^{4} (R)\rightarrow
H_{I\cap J} ^{4} (R)\rightarrow \ldots,$$ the left term is zero,
whereas the middle term is not. It follows that the right term is
non zero, too. This implies that $pd\ I(G) = cd\ I(G) \geq 4$. On
the other hand, by Lemma~\ref{lemma}, the elements
$x_1x_2,x_2x_3+x_1x_3,x_1y_2,x_1y_3+y_2y_3$ generate $I(G)$, up to
radical, which shows that $ara\ I(G) \leq 4$. It follows that $pd
\ I(G) = ara\ I(G) = 4$, which proves the claim in this case. \\
 Without loss of generality we can thus assume that $m\geq 6$. We can write $\varepsilon(G)=E_1\cup E_2$,
where
$E_1=\varepsilon(\{x_1,x_2\},\{x_3,x_4\},\ldots,\{x_m,x_1\},\{x_1,y_2\},\ldots,\{y_{n},x_1\};V)$
and
 $E_2=\varepsilon(\{x_2,x_3\};V)$. We have that $E_1\cap
E_2=\varepsilon(\{x_1\},\{x_4\},\{x_5,x_6\},\ldots,\\\{x_{m-1},x_m\},\{y_2,y_3\},\ldots,\{y_{n-1},y_n\})$.
By \cite[Lemma 3.5]{Jac2}) it follows that ${\tilde{H}}_i(E_1\cap
E_2)=\tilde{H}_{i-2}(\varepsilon(L_{m-4} \cup L_{n-1}))$, for all $i$. We also
have that $E_1=\varepsilon(H_1)$, where $H_1$ is the union of $C_n$
and the paths $x_3\ldots x_{m}x_1$ and $x_1x_2$.
Applying
Lemma~\ref{deg one} for $v=x_2$ (and $u=x_1$, so that
$N(u)=\{x_2,x_m,y_2,y_n\}$, we obtain that, for all $i$,
${\tilde{H}}_i(E_1)=\tilde{H}_{i-4}(\varepsilon(L_{m-3}\cup
L_{n-3}))$. Thus, by Theorem C, we deduce that, for all $i$,
${\tilde{H}}_i(E_1)=\tilde{H}_{i-4-\frac{2(n-3)}{3}}(\varepsilon(L_{m-3}))$.
According to Theorem D and (*) it is non zero if and only if
$i=\frac{2|V|+2}{3}-2$. Applying Theorem C, Theorem D and (*) we
also get that ${\tilde{H}}_i(E_1\cap E_2)\neq 0$ if and only if
$i=\frac{2|V|+2}{3}-4$. By the Mayer-Vietoris sequence we conclude
that ${\tilde{H}}_i(\varepsilon(G))\neq 0$ for
$i=\frac{2|V|+2}{3}-2$. In view of (*) we deduce that $pd \ I(G)\geq
\frac{2|V|+2}{3}$, as claimed.
\\

\textbf{Case $2.2$.} Let  $m\equiv 1$.
\\
We have $n\equiv 2$. Consider the induced subgraph $H'$ on
$V\setminus \{y_2\}$. Then $H'$ is the union of $C_m$ and the path
$L_{n-1}:y_3\ldots y_n y_1$. Applying Corollary~\ref{empty} to the path
$L_{n-2}: y_3\ldots y_n$, we obtain
${\tilde{H}}_{i}(\varepsilon(H'))={\tilde{H}}_{i-\frac{2(n-2)}{3}}(\varepsilon(C_m))$, for all $i$.
By Theorem B and (*) it is non zero if and only if
$i=\frac{2|V|-1}{3}-2$. In view of (*) we deduce that $pd \ I(G)\geq
\frac{2|V|-1}{3}$, as claimed.
\\

Now we find an upper bound for the arithmetical rank in each case. In the rest of the proof, we will refer to the polynomials $q_i$ introduced in Propositions~\ref{0}, \ref{1} and~\ref{2}; in each case, the polynomial $q'_i$ will be the one obtained from $q_i$ by replacing each variable $x_j$ with $y_j$.
\\
\\
In part $(a)$, for $m\equiv 2$, by Proposition~\ref{2} the sequence
$A_m:\ q_0,\ldots, q_{\frac{2(m-2)}{3}}$, generates $I(C_m)$, up to
radical and by Proposition~\ref{0} the sequence $A_n:\ q'_0,\ldots,
q'_{\frac{2n}{3}-1}$, generates $I(C_n)$, up to radical. Since
$I(G)=I(C_m)+I(C_n)$, the following sequence generates $I(G)$, up to
radical: $B:\ q_0,\ldots, q_{\frac{2(m-2)}{3}},q'_0,\ldots,
q'_{\frac{2n}{3}-1}$. This implies that $ara \ I(G)\leq
\frac{2|V|+1}{3}$.
\\
If $m\equiv 1$, then, by Proposition~\ref{1} the sequence $A_m:\
q_0,\ldots, q_{\frac{2(m-1)}{3}}$ generates $I(C_m)$, up to radical
and the sequence $A_n:\ q'_0,\ldots, q'_{\frac{2(n-1)}{3}}$
generates $I(C_n)$, up to radical. The summand
$x_1x_{m}$ of $q_1$ divides the product of the monomials
$q_{\frac{2(m-1)}{3}}=x_{m-1}x_m$ and $q'_0=y_1y_2=x_1y_2$.
Thus by Lemma~\ref{lemma} the sequence
 $B:\
q_0, q_1,  q_{\frac{2(m-1)}{3}}+q'_0, q_2, \ldots, q_{\frac{2(m-1)}{3}-1}, q'_1, \ldots, q'_{\frac{2(n-1)}{3}}$
of length $\frac{2|V|+1}{3}$ generates $I(G)$ up to radical. This
implies that $ara \ I(G)\leq \frac{2|V|+1}{3}$.
\\
\\
In part $(b)$, for $m\equiv 2$, according to Proposition~\ref{2},
the sequence $A_m:\ q_0,\ldots, q_{\frac{2(m-2)}{3}}$, generates
$I(C_m)$, up to radical and the sequence $A_n:\ q'_0,\ldots,
q'_{\frac{2(n-2)}{3}}$, generates $I(C_n)$, up to radical. The
sequence $B$ formed by the union of these two sequences generates
$I(G)$, up to radical. This implies that $ara \ I(G)\leq
\frac{2|V|}{3}$.
\\
If $m\equiv 1$, then, by Proposition~\ref{1} the sequence $A_m:\
q_0,\ldots, q_{\frac{2(m-1)}{3}}$ generates $I(C_m)$, up to radical
and by Proposition~\ref{0} the sequence $A_n:\ q'_0,\ldots,
q'_{\frac{2n}{3}-1}$ generates $I(C_n)$, up to radical. Thus by Lemma~\ref{lemma} the sequence
$B:\ q_0, q_1, q_{\frac{2(m-1)}{3}}+q'_0, q_2,\ldots, q_{\frac{2(m-1)}{3}-1}, q'_1, \ldots,
q'_{\frac{2n}{3}-1}$
 of length $\frac{2|V|}{3}$, generates $I(G)$, up to radical. So $ara \ I(G)\leq
\frac{2|V|}{3}$.
\\
\\
In part $(c)$, if $m\equiv 0$, then consider the sequence $B:\
q_0,\ldots,q_{\frac{2m}{3}-1}, q'_0,\ldots,q'_{\frac{2n}{3}-1}$,
where $A_m:\ q_0,\ldots,q_{\frac{2m}{3}-1}$ generates $I(C_m)$ and
$A_n:\ q'_0,\ldots,q'_{\frac{2n}{3}-1}$ generates $I(C_n)$, up to
radical, by Proposition~\ref{0}. This implies that $ara \ I(G)\leq
\frac{2|V|+2}{3}$.
\\
If $m\equiv 1$, then, by Proposition~\ref{1}, the sequence $A_m:\
q_0,\ldots,q_{\frac{2(m-1)}{3}}$, generates $I(C_m)$, up to radical.
By Proposition~\ref{2}, the sequence $A_n:\
q'_0,\ldots,q'_{\frac{2(n-2)}{3}}$, generates $I(C_n)$, up to
radical.  Thus by
Lemma~\ref{lemma} the sequence $B:\ q_0, q_1,
q_{\frac{2(m-1)}{3}}+q'_0,
q_2,\ldots,q_{\frac{2(m-1)}{3}-1},q'_1,\ldots, q'_{\frac{2(n-2)}{3}}
$, generates $I(G)$, up to radical. This implies that $ara \
I(G)\leq \frac{2|V|-1}{3}$. }\end{proof}
%%%%%%%%%%%%%%%%%%%%%%%%%%%%%%%%%%%%%%%%%%%%%%%%%%%%%%%%%%%%%%%%%%%%%
\begin{thm}\label{pd path}
Let $G$ be the graph formed by two cycles $C_n$ and $C_m$
with a path joining a vertex of $C_n$ to a vertex of $C_m$. Then the
following hold:
\\
\\
$(a)$ If $|V(G)|\equiv 1$, then $pd\ I(G)= ara\
I(G)=\frac{2|V(G)|-2}{3}$, whenever $m\equiv 2,\ n\equiv 2$. Otherwise,
$ pd\ I(G)= ara\ I(G)=\frac{2|V(G)|+1}{3}$.
\\
\\
$(b)$ If $|V(G)|\equiv 0$, then $pd\ I(G)= ara\
I(G)=\frac{2|V(G)|}{3}$.
\\
\\
$(c)$ If $|V(G)|\equiv 2$ and $m, n\equiv 0$ or $1$, then $ pd\
I(G)= ara\ I(G)=\frac{2|V(G)|+2}{3}$. Otherwise, $pd\ I(G)= ara\
I(G)= \frac{2|V(G)|-1}{3}$.
\end{thm}
\begin{proof}{
 Let $V=V(G)$. Consider the labeling
for $V$ such that $V(C_n)=\{y_1,y_2,\ldots,\\y_n\}$,
$V(C_m)=\{x_1,x_2,\ldots,x_m\}$ and let $P:\ z_0 z_1\ldots  z_k
z_{k+1}$ be the path in $G$, where $z_0=x_1$ and $z_{k+1}=y_1$. We
compute the reduced homology of $G$ of degree $|V|$. Up to
exchanging $m$ and $n$, we have the following cases.
\\

{\bf Case 1}. Let $k\equiv 2$.
\\
We can write
$\varepsilon(G)=\varepsilon(\{x_1,x_2\},\ldots,\{x_m,x_1\},\{y_1,y_2\},
\ldots,\{y_{n},y_1\},\{x_1,z_1\},\{z_1,z_2\},\\\ldots,\{z_{k-1},z_k\},\{z_k,y_1\};V)=E_1\cup
E_2$, where\\
$E_1=\varepsilon(\{x_1,x_2\},\ldots,\{x_m,x_1\},\{y_1,y_2\},
\ldots,\{y_{n},y_1\},\{x_1,z_1\},\{z_1,z_2\},\ldots,\{z_{k-1},z_k\};\\V)$
and $E_2=\varepsilon(\{z_k,y_1\};V)$. The intersection of these
simplicial complexes is \\$E_1\cap
E_2=\varepsilon(\{y_2\},\{y_n\},\{z_{k-1}\},
\{x_1,x_2\},\ldots,\{x_m,x_1\},\{y_3,y_4\},
\ldots,\{y_{n-2},y_{n-1}\},\\\{x_1,z_1\},\{z_1,z_2\},\ldots,\{z_{k-3},z_{k-2}\};V)$
(see \cite[Lemma 3.4]{Jac2}).
\\
By \cite[Lemma 3.5]{Jac2} it follows that
$$\tilde H_i(E_1\cap E_2)=\tilde
H_{i-3}(\varepsilon(H_1\cup L_{n-3})),$$
 for all $i$, where $H_1$ is the induced
subgraph on $V\setminus (V(C_n)\cup\{z_{k-1},z_k\})$. Applying
Corollary~\ref{empty} to the path $L_{k-2}:\ z_1\ldots z_{k-2}$, we
have
$${\tilde{H}}_{i}(E_1\cap
E_2)=\tilde{H}_{i-\frac{2(k-2)}{3}-3}(\varepsilon(C_m\cup L_{n-3}
)),$$
for all $i$. Since $E_2$ is a simplex, ${\tilde{H}}_i( E_2)=0$ for all
$i$.
\\
Applying Corollary~\ref{empty} to the path $L_{k+1}:\ z_0\ldots
z_{k}$, we have that, for all $i$,
$${\tilde{H}}_{i}(E_1)=\tilde{H}_{i-\frac{2(k+1)}{3}}(\varepsilon(
L_{m-1}\cup C_{n})).$$
\\
{\bf Case 1.1} Let $n\equiv 1$.
\\
By Theorem C, since $n-3\equiv 1$, we have that
$$\tilde H_i(E_1\cap E_2)=0$$
for all $i$.
The Mayer-Vietoris sequence then implies that $\tilde H_i(\varepsilon(G))=\tilde H_i(E_1)$ for all $i$. Moreover, in view of Proposition~\ref{1}, $I(C_n)$ is generated, up to radical, by the sequence $A_n:\
q'_0,\ldots,q'_{\frac{2(n-1)}{3}}$.
\\
\\
{\bf Case 1.1.1} Let $m\equiv 1$ or $m\equiv 0$.
\\
First suppose that $m\equiv 1$. Then $|V|\equiv 1$. From Theorem C we have that, for all $i$,
$\tilde H_i(E_1)=\tilde H_{i-\frac{2(k+1)}3-\frac{2(m-1)}3}(\varepsilon(C_n))$.
In view of Theorem B and (*) it follows that $\tilde H_i(E_1)\neq0$ if and only if
$i=\frac{2\vert V\vert+1}3-2$. Thus $\tilde H_i(\varepsilon(G))\neq0$ for
$i=\frac{2\vert V\vert+1}3-2$, which, by (*), implies that
$pd\ I(G)\geq \frac{2\vert V\vert+1}3$. \\
By Lemma~\ref{lemma}, the sequence $B:\ q_0, q_1,
q_{\frac{2(m-1)}{3}}+x_1z_1, q_2,\ldots, q_{\frac{2(m-1)}{3}-1},
z_ky_1, z_{k-1}z_k\\+q'_0, q'_1,\ldots, q'_{\frac{2(n-1)}{3}},
z_2z_3, z_1z_2+z_3z_4,\ldots,
z_{k-3}z_{k-2},z_{k-4}z_{k-3}+z_{k-2}z_{k-1}$ of length
$\frac{2|V|+1}{3}$, generates $I(G)$, up to radical, where $ A_m:\
q_0,\ldots, q_{\frac{2(m-1)}{3}}$ generates $I(C_m)$, up to radical
by Proposition~\ref{1}. Thus $ara\ I(G)\leq \frac{2|V|+1}{3}$.
\\
Now suppose that $m\equiv 0$. Then $|V|\equiv 0$. From Theorem C, Theorem B and (*) we deduce that $\tilde H_i(E_1)\neq0$ if and only if $i=\frac{2\vert V\vert}3-2$. Thus $\tilde H_i(\varepsilon(G))\neq0$ for $i=\frac{2\vert V\vert}3-2$, which, by (*), implies that $pd\ I(G)\geq \frac{2\vert V\vert}3$. \\
By Lemma~\ref{lemma}, the sequence $B:\ q'_0,
q'_1,q'_{\frac{2(n-1)}{3}}+y_1z_k, q'_2,\ldots,
q'_{\frac{2(n-1)}{3}-1},z_1x_1, z_{1}z_2+q_0, q_1,\ldots,
q_{\frac{2m}{3}-1}, z_3z_4, z_2z_3+z_4z_5,\ldots,
z_{k-2}z_{k-1},z_{k-3}z_{k-2}+z_{k-1}z_{k}$ of length
$\frac{2|V|}{3}$, generates $I(G)$, up to radical, where $ A_m:\
q_0,\ldots, q_{\frac{2m}{3}}$ generates $I(C_m)$, up to radical by
Proposition~\ref{0}. Therefore,
 we have $pd \ I(G)=ara\ I(G)= \frac{2|V|}{3}$.
\\
\\
{\bf Case 1.1.2} Let $m\equiv 2$.
\\
In this case $|V|\equiv 2$. Consider the induced subgraph $H_2$ on
$V\setminus \{z_k\}$. We have
$H_2=H'\cup C_n$, where $H'$ is the induced subgraph on $V\setminus
(V(C_n)\cup\{z_k\})$. By Theorem A
 we have that, for all $i$,
${\tilde{H}}_{i}(\varepsilon(H_2))={\tilde{H}}_{i-\frac{2n+1}{3}}(\varepsilon(H'))$.
 If we apply Corollary~\ref{empty} along the path $L_{k-2}: z_2\ldots z_{k-1}$, and then Lemma~\ref{deg one} for $v=z_1$, we further get, for all $i$,
${\tilde{H}}_{i}(\varepsilon(C_n\cup H'))={\tilde{H}}_{i-\frac{2n+1}{3}-\frac{2(k-2)}{3}-3}(\varepsilon(L_{m-3}))$,
which, by Theorem D and (*), is non zero in $i=\frac{2|V|-1}{3}-2$. By (*) this
implies $pd\ I(G)\geq \frac{2|V|-1}{3}$.
\\
By Lemma~\ref{lemma}, the sequence $B:\ z_1x_1, z_{1}z_2+q_0,
q_1,\ldots, q_{\frac{2(m-2)}{3}}, z_{3}z_{4},
z_{2}z_{3}+z_{4}z_{5},\ldots, z_{k-2}z_{k-1},
z_{k-3}z_{k-2}+z_{k-1}z_{k}, q'_0,q'_1,z_ky_1+q'_{\frac{2(n-1)}{3}},
q'_2,\ldots, q'_{\frac{2(n-1)}3-1}$ of length $\frac{2\vert
V\vert-1}3$ generates $I(G)$, up to radical, where the sequence
$A_n:\ q'_0,\ldots,q'_{\frac{2(n-1)}{3}}$ generates $I(C_n)$, up to
radical, by Proposition~\ref{1}. This shows that $ara\ I(G)\leq
\frac{2|V|-1}{3}$.
\\
\\
{\bf Case 1.2} Let $n\equiv 2$.
\\ By Theorem C, since $n-3\equiv 2$, we have that
$$\tilde H_i(E_1\cap E_2)=\tilde H_{i-\frac{2(k-2)}3-3-\frac{2(n-3)-1}3}(\varepsilon(C_m))$$
for all $i$. Moreover, by Theorem A,
$$\tilde H_i(E_1)=\tilde H_{i-\frac{2(k+1)}3-\frac{2n-1}3}(\varepsilon(L_{m-1})),$$
for all $i$.
 Moreover, in view of Proposition~\ref{2}, $I(C_n)$ is generated, up to radical, by the sequence $A_n:\
q'_0,\ldots,q'_{\frac{2(n-2)}{3}}$.
\\
\\
{\bf Case 1.2.1} Let $m\equiv 0$.
\\ In this case $|V|\equiv 1$. In view of Theorem B and (*) we have that
$\tilde H_i(E_1\cap E_2)\ne0$ if and only if $i=\frac{2\vert V\vert+1}3-3$,
in which case the homology group is $K^2$, and, according to Theorem D and (*),
we have that $\tilde H_i(E_1)\ne0$ if and only if $i=\frac{2\vert V\vert+1}3-3$,
in which case the homology group is $K$. From the Mayer-Vietoris sequence it then follows
that $\tilde H_i(\varepsilon(G))\neq0$ for $i=\frac{2\vert V\vert+1}3-2$, which, by (*),
implies that $pd\ I(G)\geq \frac{2\vert V\vert+1}3$. \\
The sequence $B:q'_0,\ldots, q'_{\frac{2(n-2)}{3}},q_0, \ldots,
q_{\frac{2m}{3}-1},z_1z_2,z_0z_1+z_2z_3,\ldots,
z_{k-1}z_{k},z_{k-2}z_{k-1}\\+z_{k}z_{k+1}$ of length $\frac{2\vert
V\vert+1}3$ generates $I(G)$, up to radical, by Lemma~\ref{lemma},
where  the sequence $A_m:\ q_0,\ldots,q_{\frac{2m}{3}-1}$ generates
$I(C_m)$, up to radical, by Proposition~\ref{0}. This implies that
$ara\ I(G)\leq \frac{2|V|+1}{3}$.
\\
\\
{\bf Case 1.2.2} Let $m\equiv 2$.
\\ In this case $|V|\equiv 0$. In view of Theorem B and (*) we have that $\tilde H_i(E_1\cap E_2)\ne0$ if and only if $i=\frac{2\vert V\vert}3-3$, and, according to Theorem C, since $m-1\equiv 1$,  we have that $\tilde H_i(E_1)=0$ for all $i$.  From the Mayer-Vietoris sequence it then follows that $\tilde H_i(\varepsilon(G))\neq0$ for $i=\frac{2\vert V\vert}3-2$, which, by (*), implies that $pd\ I(G)\geq \frac{2\vert V\vert}3$. \\
By Lemma~\ref{lemma}, the sequence $B:\ q_0,\ldots,
q_{\frac{2(m-2)}{3}},q'_0, \ldots,
q'_{\frac{2(n-2)}{3}},z_1z_2,z_0z_1+z_2z_3,\ldots,
\\z_{k-1}z_{k},z_{k-2}z_{k-1}+z_{k}z_{k+1}$ of length $\frac{2\vert V\vert}3$ generates $I(G)$, up
to radical, where the sequence $A_m:\
q_0,\ldots,q_{\frac{2(m-2)}{3}}$ generates $I(C_m)$, up to radical,
 by Proposition~\ref{2}. This implies that
$ara\ I(G)\leq\frac{2|V|}{3}$.
\\
\\
{\bf Case 1.3} Let $n\equiv m\equiv 0$.
\\ In this case $|V|\equiv 2$. In view of Theorem C, since $n-3\equiv0$, we have that, for all $i$,
$$\tilde H_i(E_1\cap E_2)=\tilde H_{i-\frac{2(k-2)}3-3-\frac{2(n-3)}3}(\varepsilon(C_m)),$$ and by  Theorem A, for all $i$,
$$\tilde H_i(E_1)=\tilde
H_{i-\frac{2(k+1)}3-\frac{2n}3}(\varepsilon(L_{m-1})).$$ According
to Theorem B and (*) it follows that
 $\tilde H_i(E_1\cap E_2)\ne0$ if and only if $i=\frac{2\vert V\vert+2}3-3$,
 in which case it is equal to $K^2$, and, in view of Theorem D and (*),
 $\tilde H_i(E_1)\ne0$ if and only if $i=\frac{2\vert V\vert+2}3-3$,
 in which case it is equal to $K$. From the  Mayer-Vietoris sequence it then follows that
 $\tilde H_i(\varepsilon(G))\neq0$ for $i=\frac{2\vert V\vert+2}3-2$, which, by (*),
 implies that $pd\ I(G)\geq \frac{2\vert V\vert+2}3$. \\
The sequence $B:\ q'_0,\ldots,q'_{\frac{2n}{3}-1}, q_0, \ldots,
q_{\frac{2m}{3}-1}, z_1z_2,x_1z_1+z_2z_3,\ldots, z_{k-1}z_{k},
z_{k-2}z_{k-1}\\+z_{k}y_1$, generates $I(G)$, up to radical, by
Lemma~\ref{lemma}, where $ A_m:\ q_0,\ldots, q_{\frac{2m}{3}-1}$
generates $I(C_m)$, up to radical, and $A_n:\ q'_0,\ldots,
q'_{\frac{2n}{3}-1}$ generates $I(C_n)$, up to radical, by
Proposition~\ref{0}. Therefore, we have $ara\ I(G)\leq
\frac{2|V|+2}{3}$.
\\
\\
{\bf Case 2} Let $k\equiv 0$.
\\
As in Case 1, we can write
$\varepsilon(G)=\varepsilon(\{x_1,x_2\},\ldots,\{x_m,x_1\},\{y_1,y_2\},
\ldots,\{y_{n},y_1\},\\\{x_1,z_1\},\{z_1,z_2\},\ldots,\{z_{k-1},z_k\},\{z_k,y_1\};V)=E_1\cup
E_2$, where\\
$E_1=\varepsilon(\{x_1,x_2\},\ldots,\{x_1,x_m\},\{y_1,y_2\},
\ldots,\{y_1,y_n\},\{x_1,z_1\},\ldots,\{z_{k-1},z_k\})$ and
\\
$E_2=\varepsilon(\{z_{k},y_1\})$.\\
 If $k=0$, then $E_1\cap E_2=\varepsilon(\{x_2\}, \{x_m\}, \{y_2\}, \{y_n\}, \{x_3,x_4\}, \ldots,
 \{x_{m-2},x_{m-1}\},\\ \{y_3,y_4\},\ldots,  \{y_{n-2},y_{n-1}\}; V\setminus\{x_1,y_1\})$, so that, by \cite[Lemma 3.5]{Jac1},
$${\tilde{H}}_i(E_1\cap E_2)={\tilde{H}}_{i-4}(\varepsilon(L_{m-3}\cup L_{n-3})),$$
for all $i$.
If $k\geq3$, then
$E_1\cap E_2=\varepsilon(\{z_{k-1}\}, \{y_2\}, \{y_n\}, \{x_1,x_2\}, \dots, \{x_m,x_1\},
 \{y_3,y_4\},\\\ldots, \{y_{n-2},y_{n-1}\}, \{x_1,z_1\}, \{z_1,z_2\}, \ldots
 \{z_{k-3},z_{k-2}\}; V\setminus\{z_k,y_1\})$, so that, by \cite[Lemma 3.5]{Jac1},
${\tilde{H}}_i(E_1\cap E_2)={\tilde{H}}_{i-3}(\varepsilon(H''\cup L_{n-3})),$
for all $i$,
 where $H''$ is the induced subgraph on $V\setminus(V(C_n)\cup \{z_{k-1},z_k\})$, i.e., it is the union of $C_m$ and the path $L_{k-1}: x_1z_1\ldots z_{k-2}$.
If we apply Corollary~\ref{empty} along the path $L_{k-3}:z_2\ldots z_{k-2}$ and then Lemma~\ref{deg one} for $v=z_1$, we deduce that, for all $i$,
$${\tilde{H}}_i(E_1\cap E_2)={\tilde{H}}_{i-6-\frac{2(k-3)}3}(\varepsilon(L_{m-3}\cup L_{n-3})),$$
which is evidently also true for $k=0$.
If $k=0$, we have that $E_1=\varepsilon(C_m\cup C_n)$, otherwise, if we apply Corollary~\ref{empty} along the path $L_k:z_1\ldots z_k$, we obtain that, for all $i$, $$\tilde H_i(E_1)=\tilde H_{i-\frac{2k}3}(C_m\cup C_n).$$
This equality is evidently also true for $k=0$. Since $E_2$
is a simplex,  we also have that ${\tilde{H}}_i( E_2)=0$ for all $i$.
\\
\\
{\bf Case 2.1} Let $n\equiv 1$.\\
In view of Theorem D (for $m=3$) and of Theorem C (for $m\geq4$),
since $n-3\equiv1$, we have that $\tilde H_i(E_1\cap E_2)=0$ for all
$i$, so that $\tilde H_i(\varepsilon(G))=\tilde H_i(E_1)$ for all
$i$. Moreover, in view of Theorem A, for all $i$,
$\tilde H_i(E_1)=\tilde H_{i-\frac{2k}3-\frac{2n+1}3}(\varepsilon(C_m)).$ \\
By Proposition~\ref{1}, the sequence $A_n: q'_0, q'_1, \ldots, q_{\frac{2(n-1)}{3}}$ generates $I(C_n)$, up to radical.\\
Let $m\equiv 1$. Then $|V|\equiv 2$. In view of Theorem B and (*),   $\tilde H_i(E_1)\ne0$ if and only if $i=\frac{2\vert V\vert+2}3-2$.  Hence $\tilde H_i(\varepsilon(G))\neq0$ for $i=\frac{2\vert V\vert+2}3-2$, which, by (*), implies that $pd\ I(G)\geq \frac{2\vert V\vert+2}3$. \\
If $k=0$, then the sequence $ q_0, q_1, q_{\frac{2(m-1)}{3}}+x_1y_1,
q_2,\ldots, q_{\frac{2(m-1)}{3}-1},q'_0,\ldots,
q'_{\frac{2(n-1)}{3}}$ of length $\frac{2|V|+2}{3}$, generates
$I(G)$, up to radical, by Lemma~\ref{lemma}, where the sequence $A_m:\
q_0,\ldots, q_{\frac{2(m-1)}{3}}$ generates $I(C_m)$, up to radical,
by Proposition~\ref{1}.
\\
If $k\geq 3$, then the sequence $B:\ q_0, q_1,
q_{\frac{2(m-1)}{3}}+x_1z_1, q_2,\ldots, q_{\frac{2(m-1)}{3}-1},
z_2z_3, z_1z_2+z_3z_4, \ldots, z_{k-1}z_k, z_{k-2}z_{k-1}+z_ky_1,
q'_0,\ldots, q'_{\frac{2(n-1)}{3}}$ of length $\frac{2|V|+2}{3}$,
generates $I(G)$, up to radical. Hence we have $ara\
I(G)\leq\frac{2|V|+2}{3}$.
\\
Let $m\equiv 2$. Then $|V|\equiv 0$. In view of Theorem B and (*),   $\tilde H_i(E_1)\ne0$ if and only if $i=\frac{2\vert V\vert}3-2$.  Hence $\tilde H_i(\varepsilon(G))\neq0$ for $i=\frac{2\vert V\vert}3-2$, which, by (*), implies that $pd\ I(G)\geq \frac{2\vert V\vert}3$. \\
\\
If $k=0$, then by Lemma~\ref{lemma}, the sequence $B:\ q'_0, q'_1,
q'_{\frac{2(n-1)}{3}}+x_1y_1, q'_2,\ldots,q'_{\frac{2(n-1)}{3}-1},
\\q_0,\ldots, q_{\frac{2(m-2)}{3}}$, of length $\frac{2|V|}{3}$,
generates $I(G)$ up to radical, where the sequence $A_m:\
q_0,\ldots, q_{\frac{2(m-2)}{3}}$ generates $I(C_m)$, up to
radical, by Proposition~\ref{2}.
\\
If $k\geq 3$, then the sequence $B:\ q'_0, q'_1,
q'_{\frac{2(n-1)}{3}}+y_1z_k,
q'_2,\ldots,q'_{\frac{2(n-1)}{3}-1},z_1z_2, x_1z_1+z_2z_3,\ldots,
z_{k-2}z_{k-1},z_{k-3}z_{k-2}+z_{k-1}z_k,  q_0,\ldots,
q_{\frac{2(m-2)}{3}}$, of length $\frac{2|V|}{3}$,  generates $I(G)$
up to radical. This implies that $ara\ I(G)\leq \frac{2|V|}{3}$.
\\
Let $m\equiv 0$. Then $|V|\equiv 1$. In view of Theorem B and (*),   $\tilde H_i(E_1)\ne0$ if and only if $i=\frac{2\vert V\vert+1}3-2$.  Hence $\tilde H_i(\varepsilon(G))\neq0$ for $i=\frac{2\vert V\vert+1}3-2$, which, by (*), implies that $pd\ I(G)\geq \frac{2\vert V\vert+1}3$. \\
\\
If $k=0$, then the sequence $B:\ q'_0, q'_1,
q'_{\frac{2(n-1)}{3}}+x_1y_1, q'_2,\ldots,
q'_{\frac{2(n-1)}{3}-1},q_0,\ldots, q_{\frac{2m}{3}-1}$ of length
$\frac{2|V|+1}{3}$, generates $I(G)$, up to radical, by
Lemma~\ref{lemma}, where the sequence $ A_m:\ q_0,\ldots,
q_{\frac{2m}{3}-1}$ generates $I(C_m)$, up to
radical, by Proposition~\ref{0}.\\
If $k\geq 3$, then the sequence $B:\ q'_0, q'_1,
q'_{\frac{2(n-1)}{3}}+y_1z_k, q'_2,\ldots, q'_{\frac{2(n-1)}{3}-1},
z_1z_2,x_1z_1+z_2z_3, \ldots,z_{k-2}z_{k-1},
z_{k-3}z_{k-2}+z_{k-1}z_k,q_0,\ldots, q_{\frac{2m}{3}-1}$ of length
$\frac{2|V|+1}{3}$, generates $I(G)$, up to radical, by
Lemma~\ref{lemma}.
 This  shows that $ ara\ I(G)\leq
\frac{2|V|+1}{3}$.
\\
\\
{\bf Case 2.2} Let $n\equiv 2$ and $m\equiv0$.
\\
In this case $|V|\equiv2$. In view of  Theorem C, Theorem D and (*), $\tilde H_i(E_1\cap E_2)\neq 0$ if and only if $i=\frac{2\vert
V\vert-1}3-2$, in which case the homology group is $K$. Moreover, in
view of Theorem A, $\tilde H_i(E_1)=\tilde
H_{i-\frac{2k}{3}-\frac{2m}{3}}(\varepsilon(C_n))$, for all $i$.
In view of Theorem B and (*),   $\tilde H_i(E_1)\ne 0$ if and only if
 $i=\frac{2\vert V\vert-1}3-2$,
in which case the homology group is $K^2$. Hence $\tilde H_i(\varepsilon(G))\neq0$ for $i=\frac{2\vert V\vert-1}3-2$, which, by (*), implies that $pd\ I(G)\geq \frac{2\vert V\vert-1}3$. \\
Let $k=0$. The sequence $B:\ x_1y_1, q_0+q'_0, q_1,\ldots,
q_{\frac{2m}{3}-1}, q'_1,\ldots, q'_{\frac{2(n-2)}{3}}$  of length
$\frac{2\vert V\vert-1}3$ generates $I(G)$, up to radical, by
Lemma~\ref{lemma}, where the sequence $A_m: q_0, q_1,\ldots, q_{\frac{2m}{3}-1}$ generates $I(C_m)$, up to radical, by Proposition~\ref{0}. \\
Let $k\geq 3$. The sequence $B:\ x_1z_1, q_0+z_1z_2, q_1,\ldots,
q_{\frac{2m}{3}-1}, z_3z_4, z_2z_3+z_4z_5, \ldots,
\\ z_{k-3}z_{k-2}, z_{k-4}z_{k-3}+z_{k-2}z_{k-1}, z_{k}y_{1}, z_{k-1}z_k+ q'_0, q'_1,\ldots,
q'_{\frac{2(n-2)}{3}}$ of length $\frac{2\vert V\vert-1}3$ generates $I(G)$, up to radical. Hence, in view of (*), we conclude that $ara\
I(G)\leq \frac{2|V|-1}{3}$.\\
\\
{\bf Case 2.3} Let $n\equiv m\equiv2$.\\
In this case $|V|\equiv 1$. Consider the induced subgraph $H_3$ on
$V\setminus \{x_2\}$. Applying Corollary~\ref{empty} to the path
$L_{m+k-2}: x_3 x_4\ldots x_1 z_1\ldots z_{k-1}$ ($L_{m-2}: x_3
x_4\ldots x_1$, if $k=0$) and Lemma~\ref{deg one} for $v=z_k$
($u=y_1$), we obtain that, for all $i$,
${\tilde{H}}_{i}(\varepsilon(H_3))={\tilde{H}}_{i-\frac{2(m+k-2)}{3}-3}(\varepsilon(L_{n-3}))$,
which, by Theorem D and (*), is non zero in $i=\frac{2|V|-2}{3}-2$.
So ${\tilde{H}}_{\frac{2|V|-2}{3}-2}(\varepsilon(G))\neq 0$ and by
(*) we have
 $pd\ I(G)\geq \frac{2|V|-2}{3}$.
\\
If $k=0$, then, by Lemma~\ref{lemma}, the sequence $B:\ x_1y_1, q_0+q'_0, q_1,\ldots,
q_{\frac{2(m-2)}{3}}, q'_1,\ldots,\\ q'_{\frac{2(n-2)}{3}}$ generates
$I(G)$, up to radical, where the sequence $A_m:\
q_0,\ldots,q_{\frac{2(m-2)}{3}}$ generates $I(C_m)$, up to radical,   and the sequence $A_n:\ q'_0,\ldots,q'_{\frac{2(n-2)}{3}}$
 generates $I(C_n)$, up to radical, by
Proposition~\ref{2}.
\\
If $k\geq 3$, then, by Lemma~\ref{lemma}, the sequence $B:\ z_ky_1, z_{k-1}z_k+q'_0,
q'_1,\ldots, q'_{\frac{2(n-2)}{3}}, z_0z_1, \\z_1z_2+q_0, q_1,\ldots,
q_{\frac{2(m-2)}{3}}, z_3z_4, z_2z_3+z_4z_5,\ldots,
z_{k-3}z_{k-2},z_{k-4}z_{k-3}+z_{k-2}z_{k-1}$ generates $I(G)$, up
to radical. Hence we have $ara\ I(G)\leq
\frac{2|V|-2}{3}$.
\\
\\
{\bf Case 2.4} Let $n\equiv m\equiv0$.\\
In this case $|V|\equiv 0$. First assume that $n=m=3$. We have that $I(G) = I \cap
J$, where $I = I(G)+(x_1y_1 z_3z_6\ldots z_k)$ and $J = (x_2,
x_3, y_2, y_3, z_1,z_2,z_4,z_5,\ldots,z_{k-5},
z_{k-4},\\ z_{k-2},z_{k-1})$.  Since $J$ is a complete intersection ideal, we have that
 $cd\ J =
4+\frac{2k}{3}$. Moreover, $I + J = (x_1y_1z_3z_6\ldots z_k, x_2,
x_3, y_2, y_3, z_1,z_2,z_4,z_5,\ldots, z_{k-5},
z_{k-4},z_{k-2},z_{k-1})$. Since $I+J$ has grade equal to
$5+\frac{2k}{3}$, by \cite[Theorem 6.2.7]{Brod} we have $H_{I+J}^{i}
(R) \neq 0$ in $i = 5+\frac{2k}{3}$ and $H_{I+J} ^{i} (R)= 0$ for
any $i<5+\frac{2k}{3}$. In the Mayer-Vietoris sequence for local
cohomology (see \cite{HT}, Section $3$)
$$\ldots\rightarrow H_{I+J} ^{4+\frac{2k}{3}}(R)\rightarrow H_{I} ^{4+\frac{2k}{3}} (R)\oplus H_{J} ^{4+\frac{2k}{3}} (R)\rightarrow
H_{I\cap J} ^{4+\frac{2k}{3}} (R)\rightarrow \ldots,$$ the left
term is zero, whereas the middle term is not. It follows that the
right term is non zero, too. This implies that $pd\ I(G) = cd\ I(G)
\geq 4+\frac{2k}{3}=\frac{2\vert V\vert}3$.
\\
So without loss of generality we may assume that $n>3$.  Then from Theorem C, since $m-3\equiv0$, we have that
${\tilde{H}}_i(E_1\cap E_2)={\tilde{H}}_{i-6-\frac{2(k-3)}3-\frac{2(m-3)}3}(\varepsilon(L_{n-3})),$ for all $i$. Hence, in view of Theorem D and (*), we have that  ${\tilde{H}}_{i}(E_1\cap
E_2)\neq 0$ only if $i=\frac{2|V|}{3}-2$, in which case this homology group is $K$.
In view of Theorem A,  Theorem B and (*) we also have that
${\tilde{H}}_{i}(E_1)\neq 0$ only if $i=\frac{2|V|}{3}-2$, in which case this homology group
 is $K^2$. The Mayer-Vietoris sequence shows
that ${\tilde{H}}_{i}(\varepsilon(G))\neq 0$ in
$i=\frac{2|V|}{3}-2$. Thus in view of (*) we deduce that $pd\
I(G)\geq \frac{2|V|}{3}$.
\\
If $k=0$, then, by Lemma~\ref{lemma}, the sequence $B:\ x_1y_1,
q_0+q'_0, q_1, \ldots, q_{\frac{2m}{3}-1},
q'_1,\ldots,\\q'_{\frac{2n}{3}-1}$ of length $\frac{2\vert V\vert}3$ generates $I(G)$, up to radical, where the sequence $A_m:\ q_0,\ldots,q_{\frac{2m}{3}-1}$
 generates $I(C_m)$, up to radical, and the sequence
 $A_n:\ q'_0,\ldots,q'_{\frac{2n}{3}-1}$
 generates $I(C_n)$, up to radical, by
Proposition~\ref{0}.
\\
If $k\geq 3$, then, by Lemma~\ref{lemma} the sequence $B:\ z_1x_1,
z_1z_2+q_0, q_1, \ldots, q_{\frac{2m}{3}-1}, z_3z_4,
\\z_2z_3+z_4z_5, \ldots,
z_{k-3}z_{k-2},z_{k-4}z_{k-3}+z_{k-2}z_{k-1}, z_ky_1,
z_{k-1}z_k+q'_0,q'_1,\ldots,q'_{\frac{2n}{3}-1}$  of length $\frac{2\vert V\vert}3$ generates $I(G)$,
up to radical.
We thus have $ara\
I(G)\leq\frac{2|V|}{3}$.
\\

{\bf Case 3} Let $k\equiv1$.\\
 We can write
$\varepsilon(G)=\varepsilon(\{x_1,x_2\},\ldots,\{x_m,x_1\},\{y_1,y_2\},
\ldots,\{y_{n},y_1\},\{x_1,z_1\},\{z_1,z_2\},\\\ldots,\{z_{k-1},z_k\},\{z_k,y_1\};V)=E_1\cup
E_2$, where\\
$E_1=\varepsilon(\{x_2,x_3\},\ldots,\{x_m,x_1\},\{y_1,y_2\},
\ldots,\{y_{n},y_1\},\{x_1,z_1\},\{z_1,z_2\},\ldots,\{z_{k-1},z_k\},\\\{z_k,y_1\};V)$
and $E_2=\varepsilon(\{x_1,x_2\};V)$. We have that $E_1=\varepsilon(L_{m+k}\cup C_n)$, where $L_{m+k}:x_2\ldots x_mx_1z_1\ldots z_k$.
\\
If $m=3$, then, by \cite[Lemma 3.4]{Jac2},
$E_1\cap
E_2=\varepsilon(\{x_3\}, \{z_1\},\{z_2,z_3\},\ldots,\{z_k,y_1\},\\\{y_1,y_2\},\ldots,\{y_{n},y_{1}\};V)$,
so that, by
\cite[Lemma 3.5]{Jac2},
$$\tilde H_i(E_1\cap E_2)=\tilde H_{i-2}(\varepsilon(H_4)),$$
for all $i$, where $H_4$ is the induced subgraph on $V\setminus (V(C_m)\cup\{z_1\})$, i.e., the union of $C_n$ and the path $L_k:z_2\ldots y_1$. If $m\geq 4$, then $E_1\cap
E_2=\varepsilon(\{x_3\},\{x_m\},\{z_1\},\{x_4,x_5\},\\\ldots,
\{x_{m-2},x_{m-1}\},\{z_2,z_3\},\ldots,\{z_k,y_1\},\{y_1,y_2\},\ldots,\{y_{n},y_{1}\};V)$, so that,  by
\cite[Lemma 3.5]{Jac2},
$$\tilde H_i(E_1\cap
E_2)=\tilde H_{i-3}(\varepsilon(L_{m-4}\cup H_4),$$ for all $i$. Since $E_2$
is a simplex, ${\tilde{H}}_i( E_2)=0$ for all $i$.
\\
\\
{\bf Case 3.1} Let $n\equiv 1$.\\
The sequence $A_n:\
q'_0,\ldots,q'_{\frac{2(n-1)}{3}}$ generates $I(C_n)$, up to
radical, by Proposition~\ref{0}.\\
\\
{\bf Case 3.1.1} Let $m\equiv 0$ or $m\equiv 2$.\\
First let $m\equiv 0$. Then $|V|\equiv 2$. If we apply
Corollary~\ref{empty} to the path $L_{m+k-1}: x_2\ldots x_mx_1\ldots
z_{k-1}$, and then Lemma~\ref{deg one} for $v=z_{k}$ we get, that, for all $i$,
${\tilde{H}}_{i}(E_1)=\tilde
H_{i-\frac{2(m+k-1)}3-3}(\varepsilon(L_{n-3})),$
which is zero for all $i$ by Theorem D and (*).
If $m\geq 6$,
applying Theorem C $(m-4\equiv2)$ and Corollary~\ref{empty} to
$E_1\cap E_2$ along the path $L_{k-1}:z_2\dots z_k$, we deduce that, for all $i$,
$${\tilde{H}}_{i}(E_1\cap E_2)=\tilde
H_{i-\frac{2(m-4)-1}3-\frac{2(k-1)}3-3}(\varepsilon(C_n)),$$ which is also true for $m=3$.\\
By Theorem B and (*),
${\tilde{H}}_{\frac{2|V|+2}{3}-3}(E_1\cap E_2)\ne0$. So by the
Mayer-Vietoris sequence
${\tilde{H}}_{\frac{2|V|+2}{3}-2}(\varepsilon (G))\ne0$ and in view of
(*) we conclude that $pd \ I(G)\geq \frac{2|V|+2}{3}$.
\\
By Lemma~\ref{lemma}, $B:\ q'_0,\ldots,q'_{\frac{2(n-1)}{3}}, x_1z_1, z_1z_2+q_0, q_1, \ldots,
q_{\frac{2m}{3}-1},  z_3z_4,
z_2z_3+z_4z_5,\\\ldots, z_{k-1}z_{k}, z_{k-2}z_{k-1}+z_{k}y_1$ of length $\frac{2\vert V\vert+2}3$
generates $I(G)$, up to radical,  where  sequence $A_m:\
q_0,\ldots,q_{\frac{2m}{3}-1}$ generates $I(C_m)$, up to radical, by Proposition~\ref{0}.
Therefore, we have $ara\ I(G)\leq \frac{2|V|+2}{3}$.\\
Now let $m\equiv 2$. In this case $|V|\equiv 1$, and $m+k\equiv 0$.
Moreover, by Theorem C, since $m-4\equiv 1$, we have that $\tilde
H_i(E_1\cap E_2)=0$ for all $i$. Hence $\tilde
H_i(\varepsilon(G))=\tilde H_i(E_1)$ for all $i$. Applying
Corollary~\ref{empty} to the path $L_{m+k}: x_2 x_3\ldots x_1
z_1\ldots z_{k}$ we obtain that, for all $i$,
${\tilde{H}}_{i}(E_1)=\tilde{H}_{i-\frac{2(m+k)}{3}}(\varepsilon(C_{n}))$.
By Theorem B and (*), ${\tilde{H}}_{i}(E_1)\neq 0$ in
$i=\frac{2|V|+1}{3}-2$. Thus by (*) we have $pd\ I(G)\geq
\frac{2|V|+1}{3}$.
\\
By Lemma~\ref{lemma}, the sequence $B:\ x_1z_1, q_0+z_1z_2,
q_1,\ldots,q_{\frac{2(m-2)}{3}},z_3z_4,z_2z_3+z_4z_5,\ldots,
z_{k-1}z_{k},z_{k-2}z_{k-1}+ z_ky_1,
q'_0,\ldots,q'_{\frac{2(n-1)}{3}}$ of length $\frac{2\vert
V\vert+1}3$ generates $I(G)$, up to radical, where  the sequence
$A_m:\ q'_0,\ldots,q'_{\frac{2(m-2)}{3}}$ generates $I(C_m)$, up to
radical, by Proposition~\ref{2}.  Thus $ara\
I(G)\leq\frac{2|V|+1}{3}$.\\
\\
{\bf Case 3.1.2} Let  $m\equiv 1$.\\
In this case $|V|\equiv 0$. Consider the induced subgraph $H_5$ on
$V\setminus \{z_k\}$. We have, for all $i$,
${\tilde{H}}_{i}(\varepsilon(H_5))={\tilde{H}}_{i}(H'''\cup C_n)$, where
$H'''$ is the induced subgraph on $V\setminus (V(C_n)\cup\{z_k\})$,
i.e., the union of $C_m$ and the path $L_k: x_1z_1\ldots z_{k-1}$.
Applying Theorem A and then Corollary~\ref{empty} to
$H'''$ along the path $L_{k-1}: z_1\dots z_{k-1}$ we have
${\tilde{H}}_{i}(\varepsilon(H_5))={\tilde{H}}_{i-\frac{2n+1}{3}-\frac{2(k-1)}{3}}(C_m),$ for all $i$, and this homology group, by Theorem B and (*), is non zero in $i=\frac{2|V|}{3}-2$. So ${\tilde{H}}_{\frac{2|V|}{3}-2}(\varepsilon(H_5))\neq 0$.
 In view of (*) we deduce that  $pd\
I(G)\geq\frac{2|V|}{3}$.
\\
The sequence $B:\ q_0, q_1, q_{\frac{2(m-1)}{3}}+z_1x_1, \ldots,q_{\frac{2(m-1)}{3}-1},
, q'_0, q'_1, q'_{\frac{2(n-1)}{3}}+z_ky_1, \ldots,\\q'_{\frac{2(n-1)}{3}-1},
z_2z_3, z_1z_2+z_3z_4, \ldots, z_{k-2}z_{k-1},
z_{k-3}z_{k-2}+z_{k-1}z_k$ of length $\frac{2|V|}{3}$, generates
$I(G)$, up to radical, by Lemma~\ref{lemma}, where the sequence $A_m:\ q_0,\ldots, q_{\frac{2(m-1)}{3}}$ generates $I(C_m)$,
up to radical, by Proposition~\ref{1}. Therefore, we have that $ara\
I(G)\leq\frac{2|V|}{3}$.
\\
\\
{\bf Case 3.2} Let $m\equiv 2$, and $n\equiv 0$ or $2$.
\\
In this case $m+k\equiv0$.
Applying Corollary~\ref{empty} to the
path $L_{m+k}: x_2 x_3\ldots x_1 z_1\ldots z_{k}$ we obtain that, for all $i$,
${\tilde{H}}_{i}(E_1)=\tilde{H}_{i-\frac{2(m+k)}{3}}(\varepsilon(C_{n}))$.
\\ Moreover, the sequence $A_m:\
q_0,\ldots,q_{\frac{2(m-2)}{3}}$ generates $I(C_m)$, up to
radical, by Proposition~\ref{2}.
\\
First let $n\equiv 0$. Then $|V|\equiv 0$ and, by Theorem B and (*),
 ${\tilde{H}}_{i}(E_1)\neq 0$
in $i=\frac{2|V|}{3}-2$. Thus by (*) we have $pd\ I(G)\geq
\frac{2|V|}{3}$.
\\
By Lemma~\ref{lemma}, the sequence $B:\ x_1z_1,
q_0+z_1z_2,q_1,\ldots,q_{\frac{2(m-2)}{3}}, z_3z_4,z_2z_3+z_4z_5,
\ldots,
z_{k-1}z_{k}, z_{k-2}z_{k-1}+z_{k}y_1, q'_0, q'_1,\ldots,
q'_{\frac{2n}{3}-1}$ of length $\frac{2\vert V\vert}3$ generates
$I(G)$, up to radical, where the sequence $A_n:\ q'_0,\ldots,
q'_{\frac{2n}{3}-1}$ generates $I(C_n)$, up to radical by
Proposition~\ref{0}. Thus $ara\ I(G)\leq\frac{2|V|}{3}$.
\\
If $n\equiv 2$, then $|V|\equiv 2$ and, by Theorem B and (*), ${\tilde{H}}_{i}(E_1)\neq 0$
in $i=\frac{2|V|-1}{3}-2$. Thus by (*) we have $pd\ I(G)\geq
\frac{2|V|-1}{3}$.
\\
By Lemma~\ref{lemma}, the sequence $B:\
q_0,\ldots,q_{\frac{2(m-2)}{3}}, z_ky_1, z_{k-1}z_k+q'_0,
q'_1,\ldots, q'_{\frac{2(n-2)}{3}}, \\z_1z_2, x_1z_1+z_2z_3, \ldots, z_{k-3}z_{k-2},
z_{k-4}z_{k-3}+z_{k-2}z_{k-1}$ of length $\frac{2\vert V\vert-1}3$
generates $I(G)$, up to radical, where the sequence $A_n:\ q'_0,\ldots, q'_{\frac{2(n-2)}{3}}$ generates $I(C_n)$,
up to radical by Proposition~\ref{2}. Thus
$ara\ I(G)\leq \frac{2|V|-1}{3}$.
\\
\\
{\bf Case 3.3} Let $n\equiv m\equiv 0$.\\
In this case $|V|\equiv1$. As in Case 1, we can write $\varepsilon(G)=E_1\cup
E_2$, where
$E_1=\varepsilon(\{x_1,x_2\},\ldots,\{x_m,x_1\},\{y_1,y_2\},
\ldots,\{y_{n},y_1\},\{x_1,z_1\},\{z_1,z_2\},\ldots,\{z_{k-1},z_k\};V)$
and $E_2=\varepsilon(\{z_k,y_1\};V)$.
\\
Applying Corollary~\ref{empty} to the path $L_{k-1} : x_1 z_1\ldots
z_{k-2}$, we have that, for all $i$, $ {\tilde{H}}_i(E_1\cap
E_2)={\tilde{H}}_{i-\frac{2(k-1)}{3}-3}(\varepsilon (L_{m-1}\cup
L_{n-3}))$. Theorem C $(m-1\equiv 2)$, Theorem D and (*) show that $
{\tilde{H}}_i(E_1\cap E_2)\neq 0$ only if $i = \frac{2|V|+1}{3}- 3$.
Since $E_2$ is a simplex, $ {\tilde{H}}_i(E_2)= 0$ for all $i$.
Applying Corollary ~\ref{empty} to $E_1$ along the path $L_{k-1} :
z_2\ldots z_k$, and once again Lemma~\ref{deg one} for $v = z_1$, we
obtain that, for all $i$, ${\tilde{H}}_i(E_1)=
{\tilde{H}}_{i-\frac{2(k-1)}{3}-3}(\varepsilon(L_{m-3}\cup C_n))$,
which by Theorem C, Theorem B and (*), is non zero only in $i =
\frac{2|V|+1}{3}- 2$. The Mayer- Vietoris sequence shows that $
{\tilde{H}}_i(\varepsilon(G))\neq 0$ in $i = \frac{2|V|+1}{3}-2$.
Thus, in view of (*), we have that $pd\ I(G)\geq
\frac{2|V|+1}{3}$.\\
By Lemma~\ref{lemma}, the sequence $B:\ x_1z_1,z_1z_2+q_0,
q_1,\ldots, q_{\frac{2m}{3}-1}, z_3z_4, z_2z_3+z_4z_5, \ldots,
\\z_{k-1}z_{k},z_{k-2}z_{k-1}+z_{k}y_1,
q'_0,\ldots,q'_{\frac{2n}{3}-1}$, generates $I(G)$, up to radical,
where the sequence $A_m:\ q_0,\ldots, q_{\frac{2m}{3}-1}$ generates
$I(C_m)$, up to radical, and the sequence $A_n:\ q'_0,\ldots,
q'_{\frac{2n}{3}-1}$ generates $I(C_n)$, up to radical, by
Proposition~\ref{0}. This implies that $pd \ I(G)=ara\ I(G)=
\frac{2|V|+1}{3}$ in this case. This completes the proof.
}\end{proof} From Theorem~\ref{pd vertex} and Theorem ~\ref{pd path}
we deduce the following corollary.
\begin{cor} Let $G$ be a bicyclic graph, then
$ara\,I(C_n)=\, pd\,I(C_n)$.
\end{cor}
%%%%%%%%%%%%%%%%%%%%%%%%%%%%%%%%%%%%%%%%%%%%%%%%%%%%%%%%%%%%%%%%%%%%%%%%%%%%%%%%%%%%%%%%

\providecommand{\bysame}{\leavevmode\hbox
to3em{\hrulefill}\thinspace}

\end{document}